\documentclass[12pt]{amsart}

\usepackage{amsmath}
\usepackage{amssymb}
\usepackage{euscript}

\newcommand{\showon}{\begin{eqnarray*}}
\newcommand{\showoff}{\end{eqnarray*}}

\newcommand{\goesto}{\rightarrow}

 \newcommand{\QQ}{\mathbb{Q}}

 \newcommand{\RR}{\mathbb{R}}

\newcommand{\Y}{\EuScript{Y}}

\begin{document}

\title{Correction to a theorem of Schoenberg}

\author{Carl Johan Ragnarsson}

\author{Wesley Wai Suen}

\author{David G. Wagner}
\address{Department of Combinatorics and Optimization\\
University of Waterloo\\
Waterloo, Ontario, Canada\ \ N2L 3G1}
\email{\texttt{wwsuen@math.uwaterloo.ca}}
\email{\texttt{dgwagner@math.uwaterloo.ca}}


\begin{abstract}
A well--known theorem of Schoenberg states that if $f(z)$ generates a PF$_r$
sequence then $1/f(-z)$ generates a PF$_r$ sequence.  We give two counterexamples
which show that this is not true, and give a correct version of the theorem.
In the infinite limit the result is sound:\  if $f(z)$ generates a PF
sequence then $1/f(-z)$ generates a PF sequence.
\end{abstract}

\maketitle

\section{The Bad News.}

Theorem $1.2$ in Chapter 8 of Karlin's book \cite{Ka} implies the following:\\

\textbf{Theorem A.}\ \emph{Let $f(z)=1+\sum_{n=1}^{\infty}a_{n}z^{n}$ 
and $g(z)=1+\sum_{n=1}^{\infty}b_{n}z^{n}$ be power series with real
coefficients such that $g(z)=1/f(-z)$.  For any positive integer $r$,
the Toeplitz matrix of $f$ is totally positive up to order $r$
if and only if the Toeplitz matrix of $g$ is totally positive up
to order $r$}.\\

Theorem A first appears in work of Schoenberg \emph{et al.} in the
early 1950s \cite{ASW,Sch1,Sch2}.  The bad news is that Theorem A is false.

First let's review the definitions.  For a power series 
$f(z)=\sum_{n=0}^{\infty}a_{n}z^{n}$, the \emph{Toeplitz matrix of 
$f$} is the infinite matrix $T[f]$, indexed by pairs of integers, with 
entries
$$T[f]_{ij}:=\left\{\begin{array}{ll}
a_{j-i} & \mathrm{if}\ j-i\geq 0,\\
0       & \mathrm{if}\ j-i<0.
\end{array}\right.$$
An infinite matrix $M$ is \emph{totally positive up to order $r$} when every minor
of $M$ of order at most $r$ is nonnegative.  This condition is abbreviated TP$_r$.
If $M$ is TP$_r$ for all $r$ then $M$ is \emph{totally positive}, abbreviated TP.

The matrix $T[f]$ is TP$_{1}$ if and only if the coefficients of 
$f(z)$ are nonnegative.  If $T[f]$ is TP$_{2}$ then the sequence of 
coefficients $a_{0}, a_{1},\ldots$ \emph{has no internal zeros:}\ if 
$0\leq h<i<j$ and $a_{h}a_{j}\neq 0$, then $a_{i}\neq 0$.
Also, if $T[f]$ is $TP_{2}$ then the sequence of coefficients
$a_{0}, a_{1},\ldots$ is \emph{logarithmically concave:}\ if $j\geq 
1$ then $a_{j}^{2}\geq a_{j-1}a_{j+1}$.  Nonegativity of the remaining
$2$--by--$2$ minors of $T[f]$ follow from these two conditions.  That 
is, the Toeplitz matrix $T[f]$ is TP$_{2}$ if and only if the sequence of 
coefficients $a_{0}, a_{1},\ldots$ is nonnegative, has no internal zeros, and is 
logarithmically concave.

Our first counterexample to Theorem A is the polynomial $f(z)=1+4z+3z^2+z^3$.
By the preceding paragraph, one sees easily that $T[f]$ is TP$_{2}$.
Elementary calculation with linear recurrence relations yields
$$g(z)=\frac{1}{1-4z+3z^{2}-z^{3}}=
1+4z+13z^2+41z^3+129z^4+406z^5+\cdots.$$
Since $129^{2}-41\cdot 406 = -5 < 0$, the Toeplitz matrix $T[g]$ is 
evidently not TP$_{2}$.  Theorem A is false.

With hindsight, one notices that the coefficients of $f(z)=1+z+2z^2$
are nonnegative, but that
$$g(z)=\frac{1}{1-z+2z^2}=1+z-z^2-3z^3-z^4+5z^5+\cdots$$
has negative coefficients.  Thus, $T[f]$ is TP$_1$ while $T[g]$ is not
TP$_1$.  This is a rather glaring counterexample to Theorem A.

\section{The Good News.}

The good news is that Theorem A can be fixed.

To do this we need a few facts about symmetric functions -- see
Macdonald \cite{Mac} for details.  Let $e_{1}, e_{2},\ldots$ and
$h_{1}, h_{2},\ldots$ be indeterminates which are algebraically 
independent over the field $\QQ$ of rational numbers, and form the
generic power series $E(t):=1+\sum_{n=1}^{\infty}e_{n}t^{n}$ and
$H(t):=1+\sum_{n=1}^{\infty}h_{n}t^{n}$.  By imposing the single 
relation $E(t)=H(-t)^{-1}$ one can determine each $e_{n}$ as a 
polyomial in the $h_{n}$-s, and conversely.  The indeterminates 
$\{h_{n}\}$ remain algebraically independent over $\QQ$, as do the 
indeterminates $\{e_{n}\}$.  The ring $\Lambda$ of polynomials with 
integer coefficients in these indeterminates is the \emph{ring of 
symmetric functions}.

Since the indeterminates $\{h_{n}\}$ are algebraically independent and
generate $\Lambda$,  a homomorphism $\varphi:\Lambda\goesto R$ from 
$\Lambda$ to another ring $R$ is determined by its values 
$\{\varphi(h_{n})\}$.  A real power series 
$f(z)=1+\sum_{n=0}^{\infty}a_{n}z^{n}$ determines such a homomorphism
$\varphi_{f}:\Lambda\goesto\RR$ by $\varphi_{f}(h_{n}):=a_{n}$.
Notice that if $g(z)=1+\sum_{n=1}^{\infty}b_{n}z^{n}$ is such that $g(z)=1/f(-z)$
then $\varphi_{f}(e_{n})=b_{n}$ and $\varphi_{g}(e_{n})=a_{n}$.

The set of all integer partitions, partially ordered by inclusion of 
Ferrers diagrams, is called \emph{Young's lattice} and denoted by $\Y$.
For $\mu\leq\lambda$ in $\Y$ there is a symmetric function 
$s_{\lambda/\mu}$ called a \emph{skew Schur function}.
Without loss of generality we may assume not only that $\mu\leq \lambda$ in $\Y$,
but also that $\mu$ has strictly fewer parts than $\lambda$ and that the largest
part of $\mu$ is strictly smaller than the largest part of $\lambda$.  We will denote
this relation by $\mu\prec\lambda$ in $\Y$.  The formulae 
we need are the Jacobi--Trudy formula and its dual form:
$$
s_{\lambda'/\mu'} = \det(e_{\lambda_{i}-i+j-\mu_{j}})
\ \ \ \
\mathrm{and}
\ \ \ \ 
s_{\lambda/\mu} = \det(h_{\lambda_{i}-i+j-\mu_{j}}).
$$
The order of these determinants is the number of parts of $\lambda$.
The notation $\lambda'$ denotes the partition conjugate to $\lambda$.
If $f(z)$ and $g(z)$ are real power series such that $g(z)=1/f(-z)$
then
$$
\varphi_f(s_{\lambda/\mu})=\varphi_g(s_{\lambda'/\mu'})
\ \ \ \
\mathrm{and}
\ \ \ \
\varphi_g(s_{\lambda/\mu})=\varphi_f(s_{\lambda'/\mu'}).
$$

Consider the submatrix $M$ of $T[f]$ supported on rows
$\{i_1<i_2<\cdots <i_r\}$ and columns $\{j_1<j_2<\cdots<j_r\}$.
If $j_{k}<i_{k}$ for any $1\leq k\leq r$ then $\det(M)=0$, so we may 
assume that $j_{k}\geq i_{k}$ for all $1\leq k\leq r$.
If $j_1=i_1$ or $j_r=i_r$ then $\det(M)$ reduces by Laplace expansion
to a smaller minor of $T[f]$.  Thus we may assume as well that $j_1>i_1$
and $j_r>i_r$.  A minor satisfying all these conditions is called 
an \emph{essential minor} of $T[f]$.  It is clear that
$T[f]$ is TP$_r$ if and only if every essential minor of $T[f]$ of
order at most $r$ is nonnegative.

Every essential minor of $T[f]$ has the form 
$\varphi_{f}(s_{\lambda/\mu})=\det(a_{\lambda_{i}-i+j-\mu_{j}})$
for some $\mu\prec \lambda$ in $\Y$.  To see this, let $\det(M)$ be 
an essential minor of $T[f]$ supported on rows
$\{i_1<i_2<\cdots <i_r\}$ and columns $\{j_1<j_2<\cdots<j_r\}$.
For each $1\leq k\leq r$ let
$\lambda_k:=j_r-i_k+k-r$.  The inequalities $\lambda_1\geq\cdots\lambda_r>0$
are easily seen, so that $\lambda$ is an integer partition with $r$ parts.
For each $1\leq k\leq r$ let $\mu_k:=j_r-j_k+k-r$.  One can check that
$\mu$ is an integer partition with at most $r-1$ parts, that $\mu\prec\lambda$
in $\Y$, and that $\det(M)=\det(a_{\lambda_{i}-i+j-\mu_{j}})$.
This construction can be reversed, so that every $\varphi_{f}(s_{\lambda/\mu})$
is an essential minor of $T[f]$.  In this way the skew Schur functions
 can be regarded as ``generic essential Toeplitz minors''.

The order of the minor $\varphi_{f}(s_{\lambda/\mu})$ of $T[f]$
is the number of parts of $\lambda$.   This implies the following:\\
\textbf{(a)}\ The Toeplitz matrix $T[f]$ is TP$_{r}$ if and only if
$\varphi_{f}(s_{\lambda/\mu})\geq 0$ for all $\mu\prec\lambda$ in 
$\Y$ for which $\lambda$ has at most $r$ parts.\\
Similarly,\\
\textbf{(b)}\ The Toeplitz matrix $T[g]$ is TP$_{r}$ if and only if
$\varphi_{g}(s_{\lambda/\mu})\geq 0$ for all $\mu\prec\lambda$ in 
$\Y$ for which $\lambda$ has at most $r$ parts.\\
If $g(z)=1/f(-z)$ then,
since $\varphi_g(s_{\lambda/\mu})=\varphi_f(s_{\lambda'/\mu'})$, condition
(b) is equivalent to\\
\textbf{(c)}\ The Toeplitz matrix $T[g]$ is TP$_{r}$ if and only if
$\varphi_{f}(s_{\lambda'/\mu'})\geq 0$ for all $\mu\prec\lambda$ in 
$\Y$ for which $\lambda$ has at most $r$ parts.\\
Or, in other words,\\
\textbf{(d)}\ The Toeplitz matrix $T[g]$ is TP$_{r}$ if and only if
$\varphi_{f}(s_{\lambda/\mu})\geq 0$ for all $\mu\prec\lambda$ in 
$\Y$ for which $\lambda$ has largest part at most $r$.\\
Comparing (a) and (d) we see that the two conditions in Theorem A are
closely related, but not equivalent.

Interpreting $\varphi_f(s_{\lambda/\mu})$ as a minor of $T[f]$,
bounding the number of parts of $\lambda$ corresponds to 
bounding the order of the minor.  What corresponds to bounding the 
largest part of $\lambda$?  For the submatrix $M$ of $T[f]$
supported on rows $\{i_1<i_2<\cdots <i_r\}$ and columns $\{j_1<j_2<\cdots<j_r\}$,
define the \emph{level} of $M$ to be $\ell:=j_r-i_1+1-r$.
The level of a minor of
$T[f]$ is the level of the submatrix of which it is the determinant.  The Toeplitz
matrix $T[f]$ is \emph{totally positive up to level $\ell$} when every
minor of $T[f]$ of level at most $\ell$ is nonnegative.  This condition is abbreviated
TP$'_\ell$.  If $T[f]$ is TP$'_\ell$ for all $\ell$ then $T[f]$ is totally
positive, TP.\\

\textbf{Theorem B.}\ \emph{Let $f(z)=1+\sum_{n=1}^{\infty}a_{n}z^{n}$ 
and $g(z)=1+\sum_{n=1}^{\infty}b_{n}z^{n}$ be power series with real
coefficients such that $g(z)=1/f(-z)$.  For any positive integer $r$,
$T[f]$ is totally positive up to \underline{level} $r$ if and only if
$T[g]$ is totally positive up to \underline{order} $r$.}\\

Notice that in the limit as $r\goesto\infty$ we get the equivalence:\
$T[f]$ is TP if and only if $T[g]$ is TP.  This is the most 
important consequence of Theorem A in the literature, and it is 
a huge relief that it survives.

\end{document}